\documentclass{amsproc}

\newtheorem{theorem}{Theorem}[section]

\newtheorem{proposition}[theorem]{Proposition}
\theoremstyle{definition}
\newtheorem{definition}[theorem]{Definition}

\theoremstyle{remark}

\numberwithin{equation}{section}

\newcommand{\Irr}[1]{\operatorname{Irr}(#1)}

\begin{document}

\title{Canonical basic sets for Hecke algebras}

\author{Nicolas Jacon}
\address{Laboratoire de Math\'ematiques Nicolas Oresme, Universit\'e de Caen, BP 5186, F 14032 Caen Cedex, France.}
\email{jacon@math.unicaen.fr}


\subjclass{Primary 20C08; Secondary  20C20}


\keywords{Modular representation theory, Hecke algebras}

\begin{abstract}
We give an explicit description of the ``canonical basic set'' for all Iwahori-Hecke algebras of finite Weyl groups in ``good'' characteristic. We obtain a complete classification of simple modules for this type of algebras.  
\end{abstract}

\maketitle

%

\section{Introduction}

Let $W$ be a finite Weyl group with set of simple reflections $S\subset{W}$ and  let $H$ be  the generic Iwahori-Hecke algebra of $W$ over $A=\mathbb{Z}[v,v^{-1}]$, where $v$ is an indeterminate. Let $u:=v^2$. $H$ has a basis $\{T_w\ |\ w\in{W}\}$ and we have the following multiplication rules. Let $s\in{S}$ and $w\in{W}$, then:
$$ T_s T_w=
    \begin{cases}
        T_{sw} & \textrm{if } l(sw)>l(w), \\
       uT_{sw}+(u-1)T_w & \textrm{if } l(sw)<l(w),\\\end{cases}$$
where $l$ is the usual length function. Let $K$ be the field of fractions of $A$ and let $\theta : A \to k$ be a homomorphism into a field $k$ such that $k$ is the field of fractions of $\theta (A)$ and such that $\theta (v)$ has finite order. 

Let $H_K:=K\otimes_A H$ and let $H_k:=k \otimes_A H$. It is known that $H_K$ is a split semi-simple algebra, isomorphic to the group algebra  $K[W]$  and that the simple $H_K$-modules are in natural bijection with the simple modules of  $K[W]$. The problem of determining a parametrization of the simple $H_k$-modules is much more complicated. To solve this problem, it is convenient to use the notion of decomposition map, which  relates the simple $H_K$-modules with the simple $H_k$-modules via a process of modular reduction. We obtain a well-defined decomposition map between the Grothendieck groups of finitely generated $H_K$-modules  and $H_k$-modules:
$$d_{\theta}:R_0(H_K)\to R_0(H_k).$$
Assume that the characteristic of $k$ is either $0$ or a good prime for $W$. Then, in \cite{G3} and in \cite{GR2}, M.Geck and R.Rouquier have defined a canonical set  $\mathcal{B}\subset{\Irr{H_K}}$ by using Lusztig's $a$-function.    
This set is called the ``canonical basic set'' and it is in natural bijection with $\Irr{H_k}$. Hence, it gives a way to parametrize the simple $H_k$-modules. Moreover, the existence of the canonical basic set implies that the decomposition matrix of $d_{\theta}$ has a lower triangular shape with $1$ along the diagonal.

In characteristic $0$, the canonical basic set has been completely described for type $A_{n-1}$ in \cite{G3}, for type $B_n$ in \cite{J3} and for type $D_n$ in \cite{G2} and \cite{J1}. Moreover,  this set  can be easily deduced for the exceptional types from the explicit tables of decomposition numbers obtained by M.Geck, K.Lux and J.M\"uller. The aim of this paper is to report these results and to show that the parametrization of $\mathcal{B}$ in characteristic $0$  holds in ``good'' characteristic. We note that the existence of ``basic sets'' has been also proved for the class of cyclotomic Hecke algebras of type $G(d,1,n)$ in  \cite{J3} and recently,  for the class of cyclotomic Hecke algebras of type $G(d,p,n)$ in \cite{GJ}.

  The paper is organized as follows. In the first part, we recall the definition of the canonical basic set.  Next, in the second part, we give some useful properties of the decomposition map and we show that the problem of determining the canonical basic set can be reduced to the case of characteristic $0$. We finally give an explicit description of the canonical basic set for all Hecke algebras of finite Weyl groups and for all specializations.
\section{Existence of canonical basic sets}
\subsection{Decomposition maps}Let $H$ be an Iwahori-Hecke algebra of a finite Weyl group $W$ over $A:=\mathbb{Z}[v,v^{-1}]$ as it is defined in the introduction. Let $K=\mathbb{Q}(v)$ and let $H_K$ be the corresponding Hecke algebra. Then $A$ is integrally closed in $K$ and  $H_K$ is a split semi-simple  algebra. Let  $\theta : A \to k$ be a specialization into a field $k$ such that $k$ is the field of fractions of $\theta (A)$ and such that $\theta (v)$ has finite order. We assume that the characteristic of $k$ is $0$ or a good prime number for $W$. Then,  there exists a discrete valuation ring $\mathcal{O}$ with maximal ideal $J(\mathcal{O})$ such that $A\subset \mathcal{O}$  and $J(\mathcal{O})\cap A=\ker{\theta}$. By \cite[Theorem  7.4.3]{GP}, we obtain a well-defined decomposition map 
$$d_{\theta}:R_0(H_K)\to R_0(H_k).$$
This is defined as follows: let $V$ be a simple $H_K$-module. Then, by \cite[section 7.4]{GP}, there exists a $H_{\mathcal{O}}$-module $\widehat{V}$ such that $K\otimes_{\mathcal{O}}\widehat{V}=V$. By reducing $\widehat V$ modulo the maximal ideal of $\mathcal{O}$, we obtain a $H_{k_{\mathcal{O}}}$-module where $k_{\mathcal{O}}$ is the residue field of $\mathcal{O}$. We obtain a map:
$$d_{\theta}':R_0(H_K)\to R_0(H_{k_{\mathcal{O}}}).$$
Since $H_k$ is split and since $k_{\mathcal{O}}$ can be seen as an extension field of $k$, we can identify $R_0 (H_{k_{\mathcal{O}}})$ with  $R_0 (H_k)$. We obtain the desired decomposition map $d_{\theta}$ between $R_0(H_K)$ and $R_0(H_k)$. Moreover, for $V\in{\Irr{H_K}}$, there exist numbers $(d_{V,M})_{M\in{\Irr{H_k}}}$ such that:
$$d_{\theta}([V])=\sum_{M\in{\Irr{H_k}}} d_{V,M}[M].$$
The matrix $(d_{V,M})_{{V\in{\Irr{H_K}}}\atop{M\in{\Irr{H_k}}}}$ is called the decomposition matrix. For more details about the construction of decomposition maps, see \cite{G4}.

\subsection{Canonical basic sets} In this part, we recall the results of \cite{G3} and \cite{GR2} which show that the above decomposition matrix has always a lower uni-triangular shape. First, we need to attach non negative integers which are called ``$a$-values'' to the simple modules of $H_K$ and $H_k$ as follows.

Let $\{C_w\}_{w\in{W}}$ be the Kazhdan-Lusztig basis of $H$. For $x,y\in{W}$, the multiplication between two elements of this basis is given by:
$$C_x C_y=\sum_{z\in{W}} h_{x,y,z} C_z$$
where $h_{x,y,z}\in{A}$ for all $z\in{W}$. For any $z\in{W}$, there is a well-defined integer $a(z)\geq 0$ such that 
\begin{align*}
& v^{a(z)}h_{x,y,z}\in{\mathbb{Z}[v]}  \textrm{ for all }x,y\in{W},\\
& v^{a(z)-1}h_{x,y,z}\notin{\mathbb{Z}[v]}  \textrm{ for some }x,y\in{W}.
\end{align*}
We obtain a function which is called the Lusztig's $a$-function:
$$\begin{array}{cccc}
   a : & W & \to & \mathbb{N}\\
     & z & \mapsto & a(z)
\end{array}
$$
Now, following \cite[Lemma 1.9]{L}, to any $M\in{\Irr{H_k}}$, we can attach an $a$-value $a(M)$ by the requirement that:
\begin{align*}
C_w.M=0 & \textrm{ for all }w\in{W}\textrm{ with } a(w)>a(M),\\
C_w.M\neq 0 & \textrm{ for some }w\in{W}\textrm{ with } a(w)=a(M).
\end{align*}
We can also attach an $a$-value $a(V)$ to any $V\in{\Irr{H_K}}$, in an analogous way.

 We can now give the theorem of existence of the canonical basic set. The main tool of the proof  is the Lusztig's asymptotic algebra.
\begin{theorem}[M.Geck \cite{G3}, M.Geck-R.Rouquier \cite{GR2}]\label{basicset} Recall that we assume that the characteristic of $k$ is $0$ or a good prime for $W$. We define the following subset of $\Irr{H_K}$:
$$\mathcal{B}:=\{V\in{\Irr{H_K}}\ |\ d_{V,M}\neq 0\ \textrm{and}\ a(V)=a(M)\ \textrm{for some }M\in{\Irr{H_k}}\}.$$
Then there exists a unique bijection 
$$\begin{array}{ccc}
 \Irr{H_k} & \to & \mathcal{B} \\
 M &  \mapsto  & V_M
\end{array}$$
such that the following two conditions hold:
\begin{enumerate}
\item For all $V_M \in{\mathcal{B}}$, we have $d_{V_M,M}=1$ and $a({V_M})=a(M)$.
\item If $V\in{\Irr{H_K}}$ and $M\in{\Irr{H_k}}$ are such that $d_{V,M}\neq 0$, then we have $a(M)\leq a(V)$, with equality only for $V=V_M$.
\end{enumerate}
The set $\mathcal{B}$ is called the canonical basic set with respect to the specialization $\theta$.
\end{theorem}
Hence, to find the elements of the canonical basic set, for each $M\in{\Irr{H_k}}$, we have to search for $V_M\in{\Irr{H_K}}$ such that $a(V_M)=a(M)$ and $d_{V_M,M}\neq 0$.

Note that a description of the set $\mathcal{B}$ would lead to a natural parametrization of the set of simple $H_k$-modules. If $H_k$ is semi-simple, we know by Tits deformation theorem that the decomposition matrix is just the identity. Hence,  we obtain the following result.
 
\begin{proposition}\label{semisimple}
Assume that $\theta$ is such that $H_k$ is a split semi-simple algebra. Then, we have:
$$\mathcal{B}=\Irr{H_K}.$$
\end{proposition}
We now want to give an explicit description of $\mathcal{B}$ in the non semi-simple case. By \cite[Theorem 7.4.7]{GP}, $H_k$ is semi-simple unless $\theta (u)$ is a root of unity. Thus, we can restrict ourselves to the case where $\theta (u)$ is a root of unity. In the next section, we will see that it is sufficient to know the canonical basic set when the characteristic of $k$ is $0$.

\section{Canonical basic sets in  positive characteristic}

In this section, we assume that $p$ is a ``good'' prime number for $W$. Let $\theta_p : A \to k_p$ be a specialization into a field $k_p$ of characteristic $p$  such that $k_p$ is the field of fractions of $\theta_p (A)$. We obtain a decomposition map
$$d_{\theta_p}:R_0 (H_K) \to R_0 (H_{k_p}).$$
Let $\mathcal{B}^p$ be the canonical basic set associated to $\theta_p$ as it is defined in Theorem \ref{basicset}.
By Proposition \ref{semisimple}, we can assume that $\theta_p (u)$ is a root of unity. We put:
$$e:=\operatorname{min}\{i\geq 2\ |\ 1+\theta_p (u)+\theta_p (u)^2+...+\theta_p (u)^{i-1}=0\}\in{\mathbb{N}}.$$

We first show that the decomposition matrix of $d_{\theta_p}$ can be obtained in two steps: one step from $u$ to a $e^{\textrm{th}}$-root of unity over $\mathbb{C}$ and another step from characteristic $0$ to characteristic $p$.

Following \cite{G5}, we denote $\mathfrak{p}:=\ker \theta_p$. Let $\Phi_e (u)\in{\mathbb{Z}[u]}$ be the $e^{\textrm{th}}$ cyclotomic polynomial. We have $\Phi_e (u)\in{\mathfrak{p}}$ and :
\begin{align*}
  \Phi_e (u) =   & \ \Phi_{2e} (v) & \textrm{if }e\textrm{ is even,} \\
    \Phi_e (u) = & \ \Phi_{e} (v)\Phi_{2e} (v) \  \textrm{and}\  \Phi_{2e} (v) =\pm \Phi_e (-v)   & \textrm{if }e\textrm{ is odd.}
\end{align*}
Thus, choosing a suitable square root of $\theta_p(u)$ in $k_{p}$, we can assume that   $\Phi_{2e} (v)\in{\mathfrak{p}}$. Let $\mathfrak{q}\subset{A}$ be the prime ideal generated by  $\Phi_{2e} (v)$. We have   $A/{\mathfrak{q}}\simeq \mathbb{Z}[\eta_{2e}]$ where $\eta_{2e}$ is a primitive $2e^{\textrm{th}}$ root of unity. Then, since $A/{\mathfrak{q}}$ is integrally closed in $k_0:=\mathbb{Q}(\eta_{2e})$ and   since $H_{k_0}$ is split, the natural map $\theta_0 : A \to A/{\mathfrak{q}}$  induces a decomposition map
$$d_{\theta_0}:R_0(H_K)\to R_0(H_{k_0}).$$
Similary, the canonical map $\pi:A/{\mathfrak{q}}\to A/{\mathfrak{p}}$ induces a decomposition map

$$d_{\pi}:R_0(H_{k_0})\to R_0(H_{k_{p}}).$$
The following result shows that   $d_{\theta_p}$ is entirely determined by $d_{\theta_0}$ and $d_{\pi}$.

\begin{proposition}[Factorization of decomposition maps, M.Geck-R.Rouquier \cite{GR1}]\label{factorization}  The following diagram is commutative:\\
\\
\begin{center}
\unitlength=1cm
\begin{picture}(6.5,2.5)
\put(2.3,2.3){\vector(1,0){2.7}} 
\put(1,2.2){$R_0(H_K)$}
\put(5,2.2){$R_0(H_{k_p})$}
\put(3,0.5){$R_0(H_{k_0})$}
\put(1.9,1.3){$d_{\theta_0}$}
\put(4.7,1.3){$d_{\pi}$}
\put(3.3,2.5){$d_{\theta_p}$}
\put(2,2.1){\vector(1,-1){1.3}}
\put(3.9,0.8){\vector(1,1){1.3}}
\end{picture}
\end{center}
If $D_{\theta_0}$,   $D_{\theta_p}$ and  $D_{\pi}$ are the decomposition matrices associated to   ${\theta_0}$,   ${\theta_p}$ and  ${\pi}$ respectively, we have:
$$D_{\theta_p}=D_{\theta_0}D_{\pi}.$$
\end{proposition}

Now, by Theorem \ref{basicset}, we have a canonical basic set $\mathcal{B}^p\subset{\Irr{H_K}}$ associated to the specialization $\theta_p$ and a canonical basic set $\mathcal{B}^0\subset{\Irr{H_K}}$ associated to    $\theta_0$.

\begin{theorem}[M.Geck-R.Rouquier \cite{GR1}, \cite{G5}]\label{number} Assume that the characteristic of $k_p$ is good. Then, we have:
$$|\Irr{H_{k_p}}|=|\Irr{H_{k_0}}|.$$
Hence, we have:
$$|\mathcal{B}^p|=|\mathcal{B}^0|$$
\end{theorem}

We can now give the main result of this section. Note that another proof of this theorem can be found in \cite[chapter 3]{J2}. The author wants to thank the referee for suggesting him this more elementary proof.

\begin{theorem}\label{restriction}  Assume that the characteristic of $k_p$ is good. Then, we have:
$$\mathcal{B}^p=\mathcal{B}^0$$
\end{theorem}
\begin{proof}
Let $\overline{M}\in{\Irr{H_{k_p}}}$ and let $V:=V_{\overline{M}}\in{\mathcal{B}^p} \subset \Irr{H_K}$ be the associated element of the canonical basic set as in Theorem \ref{basicset}. In $R_0 (H_{k_p})$, we have:
$$d_{\theta_p}([V])=[\overline{M}]+\textrm{lower terms with respect to }a-\textrm{value}$$
with $a(V)=a(\overline{M})$. Then, by the factorization of the decomposition map,  there is a simple $H_{k_0}$-module $M$ such that:
\begin{itemize}
\item   $[M]\in{R_0 (H_{k_0})}$ appears  in $d_{\theta_0}([V])$ with non zero coefficient,
\item $[\overline{M}]\in{R_0 (H_{k_p})}$ appears in $d_{\pi}([M])$ with non zero coefficient.
\end{itemize}
Then
$$a(\overline{M})=a(V)\geq a(M) \geq a(\overline{M}),$$
together with the characterization of $V_{M}$ in Theorem \ref{basicset} implies that $V=V_{M}\in{\mathcal{B}^0}$.  By Theorem \ref{number}, we have $|\mathcal{B}|=|\mathcal{B}^p|$. Thus, we obtain:
 $$\mathcal{B}=\mathcal{B}^p.$$
\end{proof}
Hence, it is sufficient to determine the canonical basic set in characteristic $0$ to describe the canonical basic set in ``good'' characteristic. 

\section{Description of the canonical basic sets}
\subsection{FLOTW multipartitions}  The canonical basic set has been explicitly computed for all types and all specializations in characteristic $0$. In this part, we recall the parametrizations of these sets and we use Theorem \ref{restriction} to compute the canonical basic set in ``good'' characteristic. 

In \cite{J3}, the existence of a ``canonical basic set'' has been proved for another type of algebras known as Ariki-Koike algebras (or cyclotomic Hecke algebras of type $G(d,1,n)$). This type of algebras can be seen as an analogue of Hecke algebras for complex reflection groups. Moreover, this canonical basic set has been explicitely described and can be parametrized by some FLOTW multipartitions arising from a crystal graph studied by Foda et al. \cite{FL} and Jimbo et al. \cite{JM}. The key of the proof is the Ariki's theory \cite{A} relating decomposition numbers with canonical basis of Fock space.

Let $v_0$, $v_1$,..., $v_{d-1}$ be integers such that $0\leq{v_0}\leq ...\leq v_{d-1}<e$. 
For $d\in{\mathbb{N}}$, we say that  $\underline{\lambda}$  is a $d$-partition of rank $n$ if:
\begin{itemize}
\item  $\underline{\lambda}=(\lambda^{(0)},..,\lambda^{(d-1)})$ where, for $i=0,...,d-1$, $\lambda^{(i)}=(\lambda^{(i)}_1,...,\lambda^{(i)}_{r_i})$ is a  partition of rank $|\lambda^{(i)}|$  such that $\lambda^{(i)}_1\geq{...}\geq{ \lambda^{(i)}_{r_i}}>0$,
\item $\displaystyle{\sum_{k=0}^{d-1}{|\lambda^{(k)}|}}=n$.
\end{itemize}
We denote by $\Pi^d_n$ the set of $d$-partitions of rank $n$. To define the FLOTW $d$-partitions, we must introduce some notations. Let  $\underline\lambda={(\lambda^{(0)},...,\lambda^{(d-1)})}$ be a $d$-partition of rank  $n$. The diagram of  $\underline{\lambda}$ is the following set:
$$[\underline{\lambda}]=\left\{ (a,b,c)\ |\ 0\leq{c}\leq{d-1},\ 1\leq{b}\leq{\lambda_a^{(c)}}\right\}.$$

The elements of this diagram are called   the nodes of  $\underline{\lambda}$. 
Let  $\gamma=(a,b,c)$ be a node of  $\underline{\lambda}$. The residue of  $\gamma$  associated to the set  $\{e;{v_0},...,{v_{d-1}}\}$ is the element of $\mathbb{Z}/e\mathbb{Z}$ defined by:
$$\textrm{res}{(\gamma)}=(b-a+v_{c})(\textrm{mod}\ e).$$

\begin{definition}
We say that  $\underline\lambda={(\lambda^{(0)} ,...,\lambda^{(d-1)})}$
is a FLOTW  $d$-partition associated to  the set  ${\{e;{v_0},...,v_{d-1}\}}$ if and only if:
\begin{enumerate}
\item for all $0\leq{j}\leq{d-2}$ and $i=1,2,...$, we have:
\begin{align*}
&\lambda_i^{(j)}\geq{\lambda^{(j+1)}_{i+v_{j+1}-v_j}},\\
&\lambda^{(d-1)}_i\geq{\lambda^{(0)}_{i+e+v_0-v_{d-1}}};
\end{align*}
\item  for all  $k>0$, among the residues appearing at the right ends of the length $k$ rows of   $\underline\lambda$, at least one element of  $\{0,1,...,e-1\}$ does not occur.
\end{enumerate}
We denote by $\Lambda^1_{\{e;{v_0},...,{v_{d-1}}\}}$ the set of FLOTW $d$-partitions associated to the set ${\{e;{v_0},...,{v_{d-1}}\}}$.
\end{definition}
Now, we are able to give the parametrizations of the canonical basic sets for all Hecke algebras of finite Weyl group and for all specializations. Let $H$ be a Hecke algebra of a finite Weyl group $W$, let $\theta:A\to k$ be a specialization into the field of fractions $k$ of $\theta (A)$ such that the characteristic of $k$ is $0$ or a good prime for $W$. Let $K$ be the field of fractions of $A$.

\subsection{Type $A_{n-1}$} Assume that $W$ is a Weyl group of type $A_{n-1}$ and that $\theta (u)$ is a primitive $e^{\textrm{th}}$-root of unity. Let $\lambda$ be a partition of rank $n$, then, we can construct an $H$-module $S^{\lambda}$, free over $A$  which is called  a Specht module (see the construction of ``dual Specht modules'' in \cite[Chapter 13]{A} in a more general setting). Moreover, we have:
$$\Irr{H_K}=\{S^{\lambda}_K\ |\ \lambda\in{\Pi^1_n}\}.$$
Now, Hecke algebras of type $A_{n-1}$ are special cases of Ariki-Koike algebras. Hence, we can use the results in \cite{J3} to find the canonical basic sets. We note that we can also find this set using results of Dipper and James \label{DJ} as it is expained in \cite[Example 3.5]{G3}. Note also that there is no bad prime number for $W=A_{n-1}$.

\begin{proposition} Assume that $W$ is a Weyl group of type $A_{n-1}$ and that $\theta (u)$ is a primitive $e^{\textrm{th}}$-root of unity. Then, we have:
$$\mathcal{B}=\{S^{\lambda}_K\ |\ \lambda\in{\Lambda^1}_{\{e;0\}},\ |\lambda|=n\}.$$
Note that:
$$\begin{array}{lll}
 \lambda\in{\Lambda^1_{\{e;0\}}} &\iff & \lambda=(\lambda_1,...,\lambda_r)\textrm{  is }e-\textrm{regular}.\\
  &\iff& \textrm{For all }i\in{\mathbb{N}},\textrm{we can't have  } \lambda_i=...=\lambda_{i+e-1}\neq 0.
\end{array}$$
\end{proposition}

\subsection{Type $B_n$} Assume that $W$ is a Weyl group of type $B_n$ and that $\theta (u)$ is a primitive $e^{\textrm{th}}$-root of unity. Let $(\lambda^{(0)},\lambda^{(1)})$ be a $2$-partition of rank $n$, then, we can construct an $H$-module $S^{(\lambda^{(0)},\lambda^{(1)})}$, free over $A$  which is called  a Specht module. Moreover, we have:
$$\Irr{H_K}=\{S^{(\lambda^{(0)},\lambda^{(1)})}_K\ |\ (\lambda^{(0)},\lambda^{(1)})\in{\Pi^2_n}\}.$$
Now, Hecke algebras of type $B_n$ are special cases of Ariki-Koike algebras. Hence, we can use the results in \cite{J3} to find the canonical basic sets. The only bad prime for type $B_n$ is $p=2$.

\begin{proposition} Assume that $W$ is a Weyl group of type $B_n$ and that $\theta (u)$ is a primitive $e^{\textrm{th}}$-root of unity. Then, we have:
\begin{itemize}
\item if $e$ is odd:
$$\mathcal{B}=\{S_K^{(\lambda^{(0)},\lambda^{(1)})}\ |\ \lambda^{(0)},\lambda^{(1)}\in{\Lambda^{1}_{\{e;0\}}},\ |\lambda^{(0)}|+|\lambda^{(1)}|=n\}.$$
\item  if $e$ is even:
$$\mathcal{B}=\{S^{(\lambda^{(0)},\lambda^{(1)})}\ |\ (\lambda^{(0)},\lambda^{(1)})\in{\Lambda^1_{\{e;1,\frac{e}{2}\}}},\ |\lambda^{(0)}|+|\lambda^{(1)}|=n\}.$$Recall that $(\lambda^{(0)},\lambda^{(1)})\in{\Lambda^1_{\{e;1,\frac{e}{2}\}}}$ if and only if:
\begin{enumerate}
\item for all  $i=1,2,...$, we have:
\begin{align*}
&\lambda_i^{(0)}\geq{\lambda^{(1)}_{i+\frac{e}{2}-1}},\\
&\lambda^{(1)}_i\geq{\lambda^{(0)}_{i+\frac{e}{2}+1}};
\end{align*}
\item  for all  $k>0$, among the residues appearing at the right ends of the length $k$ rows of   $(\lambda^{(0)},\lambda^{(1)})$, at least one element of  $\{0,1,...,e-1\}$ does not occur.
\end{enumerate}
\end{itemize}
\end{proposition}

\subsection{Type $D_n$} Assume that $W$ is a Weyl group of type $D_n$ and that $\theta (u)$ is a primitive $e^{\textrm{th}}$-root of unity. Then, $H$ can be seen as a subalgebra of an Hecke algebra $H_1$ of type $B_n$ with the following diagram (see \cite{G6}). 
\\
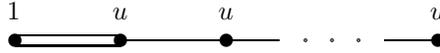
\begin{center}
\begin{picture}(240,20)
\put( 50,10){\circle*{5}}
\put( 47,18){$1$}
\put( 50,8){\line(1,0){40}}
\put( 50,12){\line(1,0){40}}
\put( 90,10){\circle*{5}}
\put( 87,18){$u$}
\put( 90,10){\line(1,0){40}}
\put(130,10){\circle*{5}}
\put(127,18){$u$}
\put(130,10){\line(1,0){20}}
\put(160,10){\circle{1}}
\put(170,10){\circle{1}}
\put(180,10){\circle{1}}
\put(190,10){\line(1,0){20}}
\put(210,10){\circle*{5}}
\put(207,18){$u$}
\end{picture}
\end{center}

The specialization $\theta$ induces a decomposition map for $H_1$:
$$d^1_{\theta}:R_0(H_{1,K})\to R_0(H_{1,k}).$$
 Similary to the equal parameter case, for all $(\lambda^{(0)},\lambda^{(1)})\in{\Pi^2_n}$,   we can construct an $H_1$-module $S^{(\lambda^{(0)},\lambda^{(1)})}$, free over $A$  which is called  a Specht module. We have:
$$\Irr{H_{1,K}}=\{S^{(\lambda^{(0)},\lambda^{(1)})}_K\ |\ (\lambda^{(0)},\lambda^{(1)})\in{\Pi^2_n}\}.$$
Now, we have an operation of restriction $\operatorname{Res}$ between the set of $H_{1,K}$-modules and the set of $H_K$-modules. For $(\lambda^{(0)},\lambda^{(1)})\in{\Pi_n^2}$, we have:
\begin{itemize}
\item if $\lambda^{(0)}\neq \lambda^{(1)}$, we have $\operatorname{Res}(S^{(\lambda^{(0)},\lambda^{(1)})}_K )\simeq \operatorname{Res}(S^{(\lambda^{(1)},\lambda^{(0)})}_K )$ and the $H_K$-module $V^{[\lambda^{(0)},\lambda^{(1)}]}:=\operatorname{Res}(S^{(\lambda^{(0)},\lambda^{(1)})}_K )$ is a simple $H_K$-module.
\item if $\lambda^{(0)}=\lambda^{(1)}$, we have $\operatorname{Res}(S^{(\lambda^{(0)},\lambda^{(1)})}_K )=V^{[\lambda^{(0)},+]}\oplus V^{[\lambda^{(0)},-]}$ where $V^{[\lambda^{(0)},+]}$ and $V^{[\lambda^{(0)},-]}$ are non isomorphic simple $H_K$-modules.
\end{itemize}
Moreover, we have:

$$\Irr{H_K}=\left\{V^{[\lambda,\mu]}\ |\ \lambda\neq{\mu},\ (\lambda,\mu)\in{\Pi^2_n}\right\}
\bigcup{\left\{V^{[\lambda,\pm]}\ |\ \lambda\in{\Pi^1_{\frac{n}{2}}} \right\}}.$$

Hecke algebras of type $B_n$ with unequal parameters are special cases of Ariki-Koike algebras. Hence, we can also define a canonical basic set for these algebras (the existence has been previously proved in \cite{G6}). Furthermore, in \cite{G6}, M.Geck has shown that the simple $H_K$-modules in  the canonical  basic set for type $D_n$ are those which appear in the restriction of the simple $H_{1,K}$-modules of the canonical basic set for type $D_n$. We obtain the following description of $\mathcal{B}$. Note that the only bad prime for type $D_n$ is $p=2$.

\begin{proposition} Assume that $W$ is a Weyl group of type $D_n$ and that $\theta (u)$ is a primitive $e^{\textrm{th}}$-root of unity. Then:
\begin{itemize}
\item if $e$ is odd, we have:
$$\begin{array}{c}
\mathcal{B}=\left\{V^{[\lambda^{(0)},\lambda^{(1)}]}\ |\ \lambda^{(0)}\neq{\lambda^{(1)}},\ \lambda^{(0)}, \lambda^{(1)}\in{\Lambda^1_{\{e;0\}}},\  |\lambda^{(0)}|+|\lambda^{(1)}|=n      \right\}\\
\bigcup{\left\{V^{[\lambda^{(0)},\pm]}\ |\ \lambda^{(0)}\in{\Lambda^1_{\{e;0\}}},\ 2|\lambda^{(0)}|=n \right\}}.\end{array}$$
\item  if $e$ is even, we have:
$$\begin{array}{c}
\mathcal{B}=\left\{V^{[\lambda^{(0)},\lambda^{(1)}]}\ |\ \lambda^{(0)}\neq{\lambda^{(1)}},\ (\lambda^{(0)},\lambda^{(1)})\in{\Lambda^1_{\{e;0,\frac{e}{2}\}}},\  |\lambda^{(0)}|+|\lambda^{(1)}|=n      \right\}\\
\bigcup{\left\{V^{[\lambda^{(0)},\pm]}\ |\ (\lambda^{(0)},\lambda^{(0)})\in{\Lambda^1_{\{e;0,\frac{e}{2}\}}},\ 2|\lambda^{(0)}|=n \right\}}.\end{array}$$
Recall that $(\lambda^{(0)},\lambda^{(1)})\in{\Lambda^1_{\{e;0,\frac{e}{2}\}}}$ if and only if:
\begin{enumerate}
\item for all  $i=1,2,...$, we have:
\begin{align*}
&\lambda_i^{(0)}\geq{\lambda^{(1)}_{i+\frac{e}{2}}},\\
&\lambda^{(1)}_i\geq{\lambda^{(0)}_{i+\frac{e}{2}}};
\end{align*}
\item  for all  $k>0$, among the residues appearing at the right ends of the length $k$ rows of   $(\lambda^{(0)},\lambda^{(1)})$, at least one element of  $\{0,1,...,e-1\}$ does not occur.
\end{enumerate}
\end{itemize}
\end{proposition}

\subsection{Exceptional types} The decomposition matrices are explicitely known for type $G_2$, $F_4$, $E_6$ and $E_7$ (see \cite{Mu}, \cite{G1},  \cite{G2} and \cite{GL}) for all specializations in characteristic $0$. Hence, in this case, it suffices to study these matrices and to use Theorem \ref{basicset}  to obtain the canonical basic in  characteristic $0$. Next,  Theorem \ref{restriction} gives the canonical basic sets in ``good'' positive characteristic. 

For type $E_8$, M\"uller has given in \cite{Mu} a set of projective $H_{\mathcal{O}}$-modules which is in bijection with the set of projective indecomposable  $H_{\mathcal{O}}$-modules (which corresponds to the columns of the decomposition matrices, see \cite[Section 7.5]{GP}). It is easy to obtain the canonical basic set by studying these   projective $H_{\mathcal{O}}$-modules.

An explicit description of the canonical basic set for all exceptional types and for all specializations can be found in \cite[Chapter 3]{J2}.

\bibliographystyle{amsalpha}

\end{document}